\documentclass[10pt,a4paper,twoside,fleqn]{article}
\usepackage{amsfonts,amssymb,amsmath,amsthm,latexsym,mathrsfs}
\usepackage{a4wide,fancyhdr}
\usepackage{graphicx,epic,eepic}
\usepackage{titlesec}
\titleformat{\section}{\large\bfseries}{\thesection}{1em}{}

\numberwithin{equation}{section}

\newtheorem{theorem}{Theorem}[section]
\newtheorem{lemma}{Lemma}[section]

\newtheorem{proposition}{Proposition}[section]

\newtheorem{corollary}{Corollary}[section]

\theoremstyle{definition}

\begin{document}

\begin{center}
{\bf\large Riesz-type inequalities and maximum flux exchange flow}
\end{center}

\begin{center}
I McGillivray\\
School of Mathematics\\
University of Bristol\\
University Walk\\
Bristol BS8 1TW\\
United Kingdom\\
e maiemg@bristol.ac.uk\\
t + 44 (0)117 3311663\\
f + 44 (0)117 9287999
\end{center}

\medskip

\begin{abstract}
\noindent Let $D$ stand for the open unit disc in $\mathbb{R}^d$ ($d\geq 1$) and $(D,\,\mathscr{B},\,m)$ for the usual Lebesgue measure space on $D$. Let $\mathscr{H}$ stand for the real Hilbert space $L^2(D,\,m)$ with standard inner product $(\cdot,\,\cdot)$. The letter $G$ signifies the Green operator for the (non-negative) Dirichlet Laplacian $-\Delta$ in $\mathscr{H}$ and $\psi$ the torsion function $G\,\chi_D$. We pose the following problem. Determine the optimisers for the shape optimisation problem
\[
\alpha_t:=\sup\Big\{(G\chi_A,\chi_A):\,A\subseteq D\text{ is open and }(\psi,\chi_A)\leq t\,\Big\}
\]
where the parameter $t$ lies in the range $0<t<(\psi,1)$. We answer this question in the one-dimensional case $d=1$. We apply this to a problem connected to maximum flux exchange flow in a vertical duct. We also show existence of optimisers for a relaxed version of the above variational problem and derive some symmetry properties of the solutions.  
\end{abstract}


\medskip


\section{Introduction}

\noindent Let $\Omega$ stand for a bounded open set in $\mathbb{R}^d$ ($d\geq 1$) and $(\Omega,\,\mathscr{B},\,m)$ for the usual Lebesgue measure space on $\Omega$. Let $\mathscr{H}$ stand for the real Hilbert space $L^2(\Omega,\,m)$ with standard inner product $(\cdot,\,\cdot)$. The letter $G$ signifies the Green operator for the (non-negative) Dirichlet Laplacian $-\Delta$ in $\mathscr{H}$ and $\psi$ the torsion function $G\,\chi_\Omega$. We pose the following problem. Determine the optimisers for the shape optimisation problem
\begin{equation}\label{modified_problem}
\alpha_t:=\sup\Big\{(G\chi_A,\chi_A):\,A\subseteq \Omega\text{ is open and }(\psi,\chi_A)\leq t\,\Big\}
\end{equation}
where the parameter $t$ lies in the range $0<t<(\psi,1)$. We show that optimisers exist for a relaxed version of this problem and derive certain symmetry properties of the solutions when $\Omega$ is replaced by the open unit ball $D$. We obtain the explicit form of the optimisers in the one-dimensional case $d=1$ for the open interval $D=(-1,\,1)$ . 

\medskip

\noindent Define
\[
V_t:=\Big\{\,\,f\in\mathscr{H}:\,0\leq f\leq 1\,m\text{-a.e. on }\Omega\text{ and }(f,\,\psi) \leq t\,\Big\}
\]
for $t$ in the range $0<t<(\psi,\,1)$ and consider the relaxed variational problem 
\begin{equation}\label{variational_problem_for_beta}
\beta_t:=\sup\Big\{\,J(f):\,f\in V_t\,\Big\},
\end{equation}
where $J(f)=(f,\,Gf)$. The first main result runs as follows.  

\begin{theorem}\label{existence_of_maximiser_for_J}
For each $t$ in the range $0<t<(\psi,\,1)$, there exists $f\in V_t$ such that $\beta_t = J(f)$.
\end{theorem}

\medskip

\noindent In case $\Omega$ is replaced by the open unit ball $D$ centred at the origin, we can say more about the symmetry properties of optimisers. In fact,

\begin{theorem}\label{cap_symmetry_of_optimisers}
Let $f\in V_t$ such that $\beta_t=J(f)$. Then $f$ possesses circular cap symmetry.
\end{theorem}

\medskip

\noindent We now turn to the one-dimensional case $d=1$ so that $D=(-1,1)$ and the torsion function $\psi$ is given explicitly by $\psi(x)=(1/2)\,(\,1-x^2\,)$ for $x\in D$. Noting that $(\psi,1)=2/3$, define $\varphi:\,D\rightarrow(0,\,2/3)$ by
$
\varphi(x):=(\chi_{(x,\,1)},\,\psi),
$
and specify $\xi_t\in(-1,\,1)$ uniquely via the relation 
\begin{equation}\label{definition_of_xi_t}
\varphi(\xi_t) = t
\end{equation}
for each $t\in(0,\,2/3)$. Set $A_t:=(\xi_t,\,1)$. Then

\begin{theorem}\label{first_main_theorem}
For any open subset $A$ in $D$ satisfying $(\psi,\,\chi_{A}) \leq t$ it holds that
\[
(G\,\chi_A,\,\chi_A) \leq (G\,\chi_{A_t},\,\chi_{A_t}),
\]
and equality occurs precisely when either $A = A_t$ or $A = -A_t$.
\end{theorem}

\medskip

\noindent The inequality is somewhat reminiscent of the Riesz rearrangement inequality: this justifies the epithet in the title. This problem has a probabilistic interpretation in so far as the function $G\,\chi_A$ is the expected occupation time in $A$ spent by absorbing Brownian motion in $D$ (associated to the Laplacian $\Delta$). The $d\geq 2$ case has not yet been resolved. It is tempting to speculate that a hyperbolic cap optimises (\ref{modified_problem}) in this case. Numerical evidence does not seem to bear this out, however \cite{Kerswell2011_1}.

\medskip

\noindent One reason why this problem is intriguing is because of its connection to maximum flux exchange flow in a vertical duct, a model of lava flow in a volcanic vent (see \cite{Kerswell2011}). In the two-dimensional case $d=2$, we imagine a configuration of two immiscible fluids in $D\times\mathbb{R}$ with different physical characteristics in a state of steady flow. The densities of the fluids are labelled $\rho$, $\rho^\prime$ and we take $\rho > \rho^\prime$. Each fluid has unit viscosity. With respect to cylindrical coordinates $(x,\,z)\in D\times\mathbb{R}$, gravity acts in the direction $(0,\,-1)$ according to the model. The pressure $p$ depends only upon $z$ and has constant gradient $\partial\,p/\partial\,z = -G$. Suppose that the fluid with density $\rho$ occupies a region in $D\times\mathbb{R}$ with cross-section $A\subseteq D$. Restricting the problem to $D$, the velocity $u$ of the components of the fluid may be described (informally) using the Navier-Stokes equation via
\[
\begin{array}{lclcl}
0 & = & \Delta\,u + G - \rho\,g & \text{ on } & A;\\
0 & = & \Delta\,u + G - \rho^\prime\,g        & \text{ on } & D\setminus A.
\end{array}
\]
Non-slip (Dirichlet) boundary conditions are imposed on the boundary of $D$. It is also assumed that $u$ and its gradient are continuous on the interface between the two regions $A$ and $D\setminus A$ (continuity of velocity and stress).  

\medskip

\noindent The parameter $G$ lies in the interval $(\rho^\prime\,g,\,\rho\,g)$. This allows the possibility of a bi-directional flow. Upon rescaling (and relabelling the velocities) we obtain the system
\begin{equation}\label{equations_for_u}
\begin{array}{lclcl}
0 & = & \Delta\,u - \lambda - 1 & \text{ on } & A;\\
0 & = & \Delta\,u - \lambda + 1 & \text{ on } & D\setminus A;
\end{array}
\end{equation}
where
\[
\lambda: = \frac{(\rho^\prime + \rho)g - 2\,G}{(\rho - \rho^\prime)\,g}\in(-1,\,1)
\]
is a proxy for the pressure gradient. Two problems arise. One is to maximise the flux $Q:=(\chi_{D\setminus A},\,u\,)$ amongst all regions $A$ which satisfy the flux balance condition $(u,\,1) = 0$ with constant $\lambda$; the other in which we optimize also over $\lambda$. In detail, we seek optimisers for the problems
\begin{eqnarray}
\label{gamma_problem}
\gamma         & := & \sup\left\{\,(\chi_{D\setminus A},\,u\,):\,(u,\,1) = 0,\,A\subseteq D\text{ open},\,\lambda\in(-1,\,1)\,\right\},\\
\label{gamma_lambda_problem}
\gamma_\lambda & := & \sup\left\{\,(\chi_{D\setminus A},\,u\,):\,(u,\,1) = 0,\,A\subseteq D\text{ open}\,\right\},
\end{eqnarray}
where in the latter $\lambda$ is fixed in the interval $(-1,\,1)$. 
It turns out that problem (\ref{modified_problem}) is closely related to the two  problems above. Note too that these problems have obvious analogues for the case $d=1$. 

\medskip

\noindent We come to our last main result. Note that the $d=2$ analogue is discussed as a marginal case in \cite{Kerswell2011}. 

\begin{theorem}\label{second_main_theorem}
In case $d=1$, 
\begin{itemize}
\item[{\it (i)}]
for each $\lambda\in(-1,\,1)$, the problem (\ref{gamma_lambda_problem}) is optimised precisely when either $A=A_{\frac{1-\lambda}{3}}$ or  $A=-A_{\frac{1-\lambda}{3}}$;
\item[{\it (ii)}]
the problem (\ref{gamma_problem}) is optimised precisely when either $A=(0,\,1)$ or $A=(-1,\,0)$ and has optimal value $1/12$. 
\end{itemize}
\end{theorem}

\medskip

\noindent We give a brief sketch of the organisation of the paper. In Section 2, we obtain existence of optimisers for the relaxed problem (\ref{variational_problem_for_beta}) and derive some symmetry properties when $\Omega$ is replaced by the ball $D$. Sections 3 to 8 deal with the proof of Theorem 1.3. Section 9 contains an application to maximum flux exchange flow (Theorem \ref{second_main_theorem}).  

\section{Existence of optimisers and symmetry in a general relaxed setting}

\noindent For the sake of clarity, we first of all remark that the (non-negative) Dirichlet Laplacian $(D(-\Delta),\,-\Delta)$ is associated with the Dirichlet form $(\mathscr{F},\,\mathscr{E})$ in $\mathscr{H}$ with form domain $\mathscr{F}:=W^{1,2}_0(\Omega)$ and 
\[
\mathscr{E}(u,\,v) = (Du,\,Dv)\hspace{1cm}(u,v\in\mathscr{F}).
\]
  
\medskip

\noindent We begin with the proof of Theorem \ref{existence_of_maximiser_for_J}.

\medskip

\noindent{\em Proof of Theorem \ref{existence_of_maximiser_for_J}.}
Let $(f_n)_{n\in\mathbb{N}}$ be a maximising sequence for $\beta_t$. 
Now, $V_t$ is weakly sequentially compact in $\mathscr{H}$. 
This follows by appeal to \cite{Lax2002} Theorem 10.2.9 due to the fact that $V_t$ is bounded, closed and convex in the reflexive Banach space $\mathscr{H}$. 
So we may assume that $(f_n)$ converges weakly in $\mathscr{H}$ to some $f\in V_t$ as $n\rightarrow\infty$ after choosing a subsequence if necessary. 

\medskip

\noindent Put $u_n:=G f_n$. Then for each $n$,
\[
\|u_n\|_{W^{1,1}_0(\Omega)}\leq \sqrt{2\,m(\Omega)}\,\|u_n\|_{W^{1,2}_0(\Omega)}.
\]
Additionally,
\[
\|u_n\|_{W^{1,2}_0(\Omega)}^2 = \mathscr{E}(u_n,\,u_n) + (u_n,\,u_n)
= (f_n,\,Gf_n) + (Gf_n,\,Gf_n)
\leq(1,\psi) + (\psi,\psi).
\]
In short, the sequence $(u_n)$ is bounded in $W^{1,1}_0(\Omega)$. In case $d\geq 2$ by the Rellich-Kondrachov compactness theorem (\cite{Evans1998} 5.7 for example), we may assume that $(u_n)$ converges in $L^1(\Omega,\,m)$ to some element $u$ after extracting a subsequence if necessary. In case $d=1$, we use Morrey's inequality (see \cite{Evans1998} 5.6.2, for example) and the Arzela-Ascoli compactness criterion to extract a uniformly convergent subsequence. The details are described in the proof of Theorem \ref{existence_of_maximiser_for_J_in_U_m_t}.  

\medskip

\noindent For each $n\in\mathbb{N}$,
\[
(u_n,\,\varphi) = (G\,f_n,\,\varphi) = (f_n,\,G\,\varphi)\text{ for all }\varphi\in\mathscr{H},
\]
which yields
\[
(u,\,\varphi) = (f,\,G\,\varphi) = (G\,f,\,\varphi)\text{ for all }\varphi\in \mathscr{H}
\]
upon taking limits. Therefore, $u=Gf$ $m$-a.e. on $\Omega$. Moreover,
\[
J(f) - J(f_n) = (\,u,\,f) - (\,u_n,\,f_n\,) = (\,u,\,f-f_n\,) + (\,f_n,\,u-u_n\,)
\]
and the right-hand side converges to zero as $n\rightarrow\infty$ in virtue of the weak respectively strong $L^1(\Omega,\,m)$ (or uniform in the case $d=1$) convergence of the sequences $(f_n)$ respectively $(u_n)$. As
$
\beta_t=\lim_{n\rightarrow\infty}J(f_n)
$
it follows that $\beta_t=J(f)$. 
\medskip
\qed

\medskip

\noindent In the remainder of this section, we replace $\Omega$ by the open unit ball $D$ in $\mathbb{R}^d$ centred at the origin. We first discuss the operation of polarisation for integrable functions on $D$ (see \cite{BrockSolynin1999} and references therein). For $\nu\in S^{d-1}$ the closed half-space $H=H_\nu$ is defined by
\[
H_\nu:=\left\{x\in\mathbb{R}^d:\,x\cdot\nu\geq 0\right\}
\]
with an associated reflection
\[
\tau_H:\,\mathbb{R}^d\rightarrow\mathbb{R}^d;\,x\mapsto x - 2\,(x\cdot\nu)\,\nu. 
\]
Refer to the collection of all these closed half-spaces by $\mathcal{H}$. The polarisation $f_H$ of $f\in L^1_+(D,\,m)$ with respect to $H\in\mathcal{H}$ is defined as follows. Choose an $m$-version of $f$, which we again denote by $f$. Set
\[
f_H(x):=
\left\{
\begin{array}{lcl}
f(x)\wedge f(\tau_H x) & \text{ for } & x\in D\cap H,\\
f(x)\vee   f(\tau_H x) & \text{ for } & x\in D\setminus H.
\end{array}
\right.
\]   
Its $m$-equivalence class is the polarisation of $f$. The definition is well-posed.  

\medskip

\noindent  The Green kernel $G(x,\,y)$ is given by
\[
G(x,\,y)=\Phi(y-x) - \Phi(|x|(y-x^*))\text{ for }(x,\,y)\in D\times D\setminus\mathtt{d},
\]
where $\Phi$ is the fundamental solution of Laplace's equation in $\mathbb{R}^d$, $\mathtt{d}$ stands for the diagonal in $D\times D$ and the decoration $*$ refers to inversion in the unit sphere. We note the inequality
\begin{equation}\label{Green_kernel_inequality}
G(x,\,y) > G(x,\,\tau_H y)\text{ for any }x,\,y\in D\cap \text{int }H,
\end{equation}
which follows from the strong maximum principle. 

\begin{theorem}\label{polarisation_and_the_functional_J}
Let $f\in L^1_+(D,\,m)$ and $H\in\mathcal{H}$. Then $J(f)\leq J(f_H)$ with equality if and only if either $f=f_H$ or $f\circ\tau_H=f_H$ $m$-a.e. on $D$.  
\end{theorem}

\noindent{\em Proof.} 
We work with an $m$-version of $f$, again denoted $f$. Define
\[
A^+:=\left\{x\in D\cap H:\,f(x)<f(\tau_H x)\right\}
\]
and similarly $B^+$ but with the strict inequality replaced by the sign $>$. Put $A^-:=\tau_H A^+$ and $A:=A^+\cup A^-$. Set $S:=D\setminus A$. In this notation,
\[
f_H=\chi_A\,f\circ\tau_H + \chi_S\,f.
\]
As a consequence,
\[
J(f_H) 
= J(\chi_A\,f\circ\tau_H) + 2\,(\chi_A\,f\circ\tau_H,\,G\chi_S f) + J(\chi_S\,f)
= J(\chi_A\,f) + 2\,(\chi_A\,f\circ\tau_H,\,G\chi_S f) + J(\chi_S\,f)
\]
and a similar identity holds for $J(f)$ but without composition with reflection. We may then write that
\[
J(f_H)-J(f)
=
2\int_{A^+}\int_{B^+}(f(\tau_H x)-f(x))(g(x,\,y)-g(\tau_H x,\,y))(f(y)-f(\tau_H y))\,m(dy)\,m(dx).
\]
It is clear from this representation with the help of (\ref{Green_kernel_inequality}) that $J(f)\leq J(f_H)$.

\medskip

\noindent In the case of equality, it holds that either $m(A^+)=0$ or $m(B^+)=0$. In the former case, $f=f_H$ while in the latter, $f\circ\tau_H = f_H$ $m$-a.e. on $D$. 

\medskip

\noindent The spherical cap symmetrisation (see \cite{vanScaftingen2004}, \cite{vanScaftingen2005}, \cite{SmetsWillem2003} for example) of $A\in\mathscr{B}$ with respect to the direction $\omega\in S^{d-1}$ is the set $A^*\in\mathscr{B}$ specified uniquely by the conditions
\[
\begin{array}{lcll}
A^*\cap\left\{0\right\}  & = & A\cap\left\{0\right\},                 & \\
A^*\cap\partial B(0,\,r) & = & B(r\omega,\,\rho)\cap\partial B(0,\,r) & \text{ for some }\rho\geq 0,\\
\sigma_r(A^*\cap\partial B(0,\,r)) & = &  \sigma_r(B(r\omega,\,\rho)\cap\partial B(0,\,r)),
\end{array}
\] 
for each $r\in(0,\,1)$. Here, $\sigma_r$ stands for the surface area measure on $\partial B(0,\,r)$. The spherical cap symmetrisation of $f\in L^1_+(D,\,m)$ (denoted $f^*$ for brevity) is defined as follows. Choose an $m$-version of $f$, which we again denote by $f$. Let $f^*$ be the unique function such that
\[
\left\{f^*>t\right\} = \left\{f>t\right\}^*\text{ for each }t\in\mathbb{R}.
\]
Its $m$-equivalence class is the polarisation of $f$. The definition is again well-posed. We also write $f^*$ as $C_\omega f$. 

\medskip

\noindent Before proving Theorem \ref{cap_symmetry_of_optimisers}, we prepare a number of lemmas. We first discuss a useful two-point inequality. We introduce the notation
\[
\begin{array}{lcl}
Q & := & \left\{(x_1,\,x_2)\in\mathbb{R}^2:\,x_1\geq 0\text{ and }x_2\geq 0\right\},\\
R & := & \left\{(x_1,\,x_2)\in Q:\,0\leq x_2 < x_1\right\},\\
S & := & \left\{(x_1,\,x_2)\in Q:\,0\leq x_1 < x_2\right\}.
\end{array}
\]
Equip $Q$ with the $\ell^1$-norm $\|x\|_1:=|x_1|+|x_2|$ where $x=(x_1,\,x_2)\in Q$. Define a mapping $\varphi:\,Q\rightarrow Q$ via
\[
(x_1,\,x_2)\mapsto(x_1\vee x_2,\,x_1\wedge x_2).
\] 
A geometric argument establishes the following lemma.

\begin{lemma}\label{two_point_inequality}
For any $x,\,y\in Q$,
\[
\|\varphi x - \varphi y\|_1\leq\|x-y\|_1
\]
with strict inequality if and only if $x\in R$ and $y\in\overline{S}$ or $x\in\overline{R}$ and $y\in S$ or the same with the r\^oles of $x$ and $y$ interchanged. 
\end{lemma}

\medskip

\noindent For $\omega\in S^{d-1}$ introduce the collection of closed half-spaces
\[
\mathcal{H}_\omega:=\left\{x\in\mathbb{R}^d:\,x\cdot\nu\geq 0\right\}.
\]

\begin{lemma}\label{inequality_involving_polarisation_and_symmetrisation}
Let $f\in L^1_+(D,\,m)$ and $\omega\in S^{d-1}$. For any $H\in\mathcal{H}_\omega$,
\begin{equation}\label{inequality_involving_polarisation_and_symmetrisation_as_equation}
\|f_H - C_\omega f\|_{L^1(D,\,m)}\leq \|f - C_\omega f\|_{L^1(D,\,m)}
\end{equation}
with strict inequality if
\[
m(\left\{f\circ \tau_H > f\right\})>0. 
\]
\end{lemma}

\noindent{\em Proof.}
Select an $m$-version of $f$, again denoted $f$. Note that $f_H^*=f^*$. By the two-point inequality Lemma \ref{two_point_inequality},
\begin{equation}\label{polarisation_and_symmetrisation_inequality}
|f_H(x)-f^*(x)| + |f_H(\tau_H x) - f^*(\tau_H x)|
\leq 
|f(x)-f^*(x)| + |f(\tau_H x) - f^*(\tau_H x)|
\end{equation}
for $x\in D\cap H$. It only remains to integrate over $D\cap H$ to obtain the inequality. 

\medskip

\noindent For each $x\in D\cap H$ the pair $(f^*(x),\,f^*(\tau_H x))$ belongs to $\overline{R}$. By Lemma \ref{two_point_inequality} the condition $(f^*(x),\,f^*(\tau_H x))\in S$ guarantees strict inequality in (\ref{polarisation_and_symmetrisation_inequality}). This observation leads to the criterion in the Lemma.  
\qed

\medskip

\noindent The next lemma is a spherical cap symmetrisation counterpart to \cite{BrockSolynin1999} Lemma 6.3, and extends \cite{vanScaftingen2004} Lemma 3.9. 

\begin{lemma}\label{polarisation_improves_distance_to_cap_symmetrisation}
Let $f\in L^1_+(D,\,m)$ and $\omega\in S^{d-1}$ and assume that $f\neq C_\omega f$. Then there exists $H\in\mathcal{H}_\omega$ such that
\[
\| f_H - C_\omega f\|_{L^1(D,\,m)}
<
\| f - C_\omega f\|_{L^1(D,\,m)}.
\]
\end{lemma}

\medskip

\noindent{\em Proof.}
For shortness, write $f^*$ for $C_\omega f$. As $f\neq f^*$ there exists $t>0$ such that
\[
m(\left\{f>t\right\}\Delta\left\{f^*>t\right\})>0.
\]
It follows that the sets $A:=\left\{f\leq t<f^*\right\}$ and $B:=\left\{f^*\leq t<f\right\}$ are disjoint and have identical positive $m$-measure. 

\medskip

\noindent We claim that there exists $H\in\mathcal{H}_\omega$ such that $m(A\cap\tau_H B)>0$. Taking this as read, on $A\cap\tau_H B$ we have that $f^*>t\geq f^*\circ\tau_H$ so that $A\cap\tau_H B\subseteq H$. Also, $f\leq t<f\circ\tau_H$ there. In short, $A\cap\tau_H B\subseteq\left\{ f\circ\tau_H>f\right\}\cap H$. So $m(\left\{ f\circ\tau_H>f\right\})>0$ and there is strict inequality in (\ref{inequality_involving_polarisation_and_symmetrisation_as_equation}) by Lemma \ref{inequality_involving_polarisation_and_symmetrisation}.  

\medskip

\noindent To prove the claim, we assume for a contradiction that $m(A\cap\tau_H B)=0$ for all $H\in\mathcal{H}_\omega$. Let $F$ be a countable dense subset in $S^{d-1}\cap H_\omega$. Then
\[
m(A\cap\bigcup_{\nu\in F}\tau_{H_\nu}B)=0. 
\]
Therefore, for all $r\in(0,\,1)$, it holds that
\[
\sigma_r(A_r\cap\tau_{H_\nu}B_r)=0\text{ for every }\nu\in F,
\]
except on a $\lambda$-null set $N$. Here, $\lambda$ stands for Lebesgue measure on the Borel sets in $\mathbb{R}$, and $A_r:=A\cap\partial B(0,\,r)$ for the section of $A$ (likewise for $B_r$). Let $\nu\in S^{d-1}\cap H_\omega$ with corresponding reflection $\tau=\tau_{H_\nu}$. Select a sequence $(\nu_j)$ in $F$ which converges to $\nu$ in $S^{d-1}$. Write $\tau_j$ for the reflection associated to closed half-space $H_{\nu_j}$. For $r\in(0,\,1)\setminus N$,
\[
|\sigma_r(A_r\cap\tau B_r)-\sigma_r(A_r\cap\tau_j B_r)|
\leq
\|\chi_B -\chi_B\circ\tau\circ\tau_j\|_{L^1(\partial B(0,\,r),\sigma_r)},
\]
and this latter converges to zero as $j\rightarrow\infty$. This is due to the fact that the special orthogonal group $SO(d)$ acts continuously on $L^1(S^{d-1},\,\sigma)$. We derive therefore that
\begin{equation}\label{vanishing_of_sections}
\sigma_r(A_r\cap\tau_{H_\nu}B_r)=0\text{ for every }\nu\in S^{d-1}\cap H_\omega
\end{equation}
for all $r\in(0,\,1)\setminus N$. 

\medskip

\noindent To conclude the argument, choose $r\in(0,\,1)\setminus N$ such that $\sigma_r(A_r)=\sigma_r(B_r)>0$. Use Lebesgue's density theorem to select a density point $x$ for $A_r$ lying in $A_r$, and choose $y$ in $B_r$ similarly. Then $f^*(x)>t\geq f^*(y)$. So there exists $\nu\in S^{d-1}\cap H_\omega$ such that with $\tau=\tau_{H_\nu}$ we have that $\tau y = x$. But this means that
\[
\lim_{\varepsilon\downarrow 0}\frac{\sigma_r(A_r\cap\tau B_r\cap B(x,\,\varepsilon))}{\sigma_r(\partial B(0,\,r)\cap B(x,\,\varepsilon))}=1,
\]
so that, in fact, $\sigma_r(A_r\cap\tau B_r)>0$, contradicting (\ref{vanishing_of_sections}). 
\qed

\medskip

\noindent{\em Proof of Theorem \ref{cap_symmetry_of_optimisers}}.
Assume for a contradiction that $f\neq C_\omega f$ for each $\omega\in S^{d-1}$. Then there exists $\omega\in S^{d-1}$ such that
\[
\delta:=\inf_{\nu\in S^{d-1}}\|f-C_\nu f\|_{L^1(D,\,m)}=
\|f-C_\omega f\|_{L^1(D,\,m)}>0.
\]
By Lemma \ref{polarisation_improves_distance_to_cap_symmetrisation} there exists $H\in\mathcal{H}_\omega$ such that
\[
\|f_H-C_\omega f\|_{L^1(D,\,m)}<\|f-C_\omega f\|_{L^1(D,\,m)}.
\]
It is plain that $f\neq f_H$. But also $f\circ\tau_H\neq f_H$, for otherwise,
\[
\|f-C_{\sigma\omega} f\|_{L^1(D,\,m)}
=\|f_H-C_{\omega} f\|_{L^1(D,\,m)}
<\|f-C_{\omega} f\|_{L^1(D,\,m)},
\]
contradicting optimality of $\omega$. It follows by Theorem \ref{polarisation_and_the_functional_J} that $J(f)<J(f_H)$ and this contradicts the optimality of $f$ in the expression for $\beta_t$. 
\qed
\section{Preliminaries for the one-dimensional problem}

\noindent In the remainder of the article we work in the one-dimensional setting where $D=(-1,\,1)$. In this context, the corresponding Green operator $G$ has kernel given by
\begin{equation}\label{formula_for_green_function}
G(x,y) =
\left\{
\begin{array}{lcl}
\frac{1}{2}(1-y)(1+x) & \text{ for } & x\leq y,\\
\frac{1}{2}(1+y)(1-x) & \text{ for } & x>y,\\
\end{array}
\right.
\end{equation}
for $x,y\in D$. We record the useful inequality
\begin{equation}\label{variation_estimate_for_Green_function}
\big|\,G(x,\,y) - G(x,\,x)\,\big|\leq\big|\,y-x\,\big|
\text{ for all }x,\,y\in D, 
\end{equation}
for future use. As noted above, the torsion function $\psi:=G\,\chi_D$ is given explicitly by $\psi(x)=(1/2)\,(\,1-x^2\,)$ for $x\in D$, and
\begin{equation}\label{inner_product_of_psi_with_1}
(1,\,\psi) = 2/3. 
\end{equation}
The Green kernel may be bounded in terms of $\psi$; that is, 
\begin{equation}\label{bound_for_G_in_terms_of_psi}
G(x,\,y)\leq\psi(x)\text{ for all }y\in D,
\end{equation}
with fixed $x\in D$. 

\medskip

\noindent For $t\in(0,\,2/3)$ introduce the shape space 
\[
U_t:=\Big\{\,\,f=\chi_A:\,A\subseteq D\text{ is open and }(f,\,\psi) \leq t\,\Big\}.
\]
We may then write
\begin{equation}\label{variational_problem_for_alpha}
\alpha_t=\sup\Big\{\,J(f):\,f\in U_t\,\Big\}.
\end{equation}
For each $t\in(0,\,2/3)$ and $m\in\mathbb{N}$ define $U^{(m)}_t$ to be the collection of all functions of the form $f=\chi_A$ where $A$ is a union of at most $m$ disjoint open intervals in $D$ with the additional requirement that $(f,\,\psi) \leq t$. We occassionally refer to the condition
\begin{equation}\label{condition_on_A}
\text{int}\,\overline{A} = A.
\end{equation}
We also introduce the variational problem
\begin{equation}\label{variational_problem_for_beta_m_t}
\alpha^{(m)}_t:=\sup\Big\{\,J(f):\,f\in U^{(m)}_t\,\Big\}.
\end{equation}
We now derive the crucial property that (\ref{variational_problem_for_beta_m_t}) attains its optimum.
 
\begin{theorem}\label{existence_of_maximiser_for_J_in_U_m_t}
For each $t\in(0,\,2/3)$ and $m\in\mathbb{N}$ there exists $f\in U^{(m)}_t$ with $(f,\,\psi) = t$ such that $\alpha^{(m)}_t = J(f)$.
\end{theorem}

\noindent{\em Proof.} Let $(f_n)_{n\in\mathbb{N}}$ be a maximising sequence for $\alpha^{(m)}_t$. Each $f_n$ may be written in the form $f_n = \sum_{j=1}^{k_n}\chi_{A_{nj}}$ for some $1\leq k_n\leq m$ where $A_{nj}=(a_{nj},\,b_{nj})$ and
\[
-1\leq a_{n1} < b_{n1} \leq a_{n2} < b_{n2} \leq \cdots \leq a_{nk_n} < b_{nk_n} \leq 1.
\]
After selecting a subsequence if necessary we may suppose that $k_n$ takes a fixed value $k$ for some $k$ between $1$ and $m$. On appeal to the Bolzano-Weierstrass theorem, we may assume (perhaps after discarding a subsequence) that $a_{nj}\rightarrow a_j$ and $b_{nj}\rightarrow b_j$ as $n\rightarrow\infty$ where
\begin{equation}\label{condition_on_as_and_bs}
-1\leq a_1 \leq b_1 \leq a_2 \leq b_2 \leq \cdots \leq a_k \leq b_k \leq 1.
\end{equation}
Set $f:=\sum_{j=1}^{k}\chi_{A_j}$ where $A_j=(a_j,\,b_j)$. By the dominated convergence theorem, $(f_n)$ converges weakly to $f$ in $\mathscr{H}$. 

\medskip

\noindent Put $u_n:=G f_n$ as before. Then the sequence $(u_n)$ is bounded in $W^{1,2}_0(D)$ as in the proof of Theorem \ref{existence_of_maximiser_for_J}. By Morrey's inequality (see \cite{Evans1998} 5.6.2, for example),
\[
\|\,u_n\,\|_{C^{0,1/2}(\overline{D})}\leq c\text{ for all }n\in\mathbb{N}
\]
for some finite constant $c$; in particular,
\[
|\,u_n(x) - u_n(y)\,| \leq c\,|\,x-y\,|^{1/2}\text{ for any }x,y\in \overline{D}
\]
and any $n\in\mathbb{N}$. Thus, $(u_n)$ forms a bounded and equicontinuous sequence in $C(\overline{D})$. By the Arzela-Ascoli compactness criterion, we may assume that $(u_n)$ converges uniformly to some $u\in C(\overline{D})$ as $n\rightarrow\infty$ after extracting a subsequence if necessary. Now continue the argument as in the proof of Theorem \ref{existence_of_maximiser_for_J} to conclude that $\alpha^{(m)}_t=J(f)$.

\medskip

\noindent We now show that $(f,\,\psi) = t$. First note that $(f,\,\psi) \leq  t$; this flows from the fact that $f$ is a weak limit of elements in $U^{(m)}_t$. Suppose for a contradiction that $(f,\,\psi) < t$. As $(f,\,\psi) < 2/3$, in (\ref{condition_on_as_and_bs}) there must exist $j=0,\ldots,k$ such that $b_j < a_{j+1}$ with the understanding that $b_0:=-1$ and $a_{k+1}:=1$. By choosing $B$ to be a suitable (semi-)open interval in $[b_j,\,a_{j+1}]$ we can arrange that the function $f_1:=f+\chi_B$ satisfies the requirement $(f_1,\,\psi)\leq t$ as well as $J(f)<J(f_1)$. This contradicts the optimality of $f$. 
\qed

\medskip

\noindent We now revisit the operation of polarisation in the one-dimensional setting. We use the letter $P$ to signify the polarisation operator with respect to the closed half-space $[0,\,\infty)$. Thus, for $f\in U_t$, the polarisation is defined by 
\begin{equation}\label{polarisation}
Pf(x) := 
\left\{
\begin{array}{lcr}
f(x)\vee f(-x)   & \text{ if } & 0 \leq x < 1,\\
f(x)\wedge f(-x) & \text{ if } & -1 < x < 0.
\end{array}
\right.
\end{equation}
Alternatively, suppose that $f=\chi_A$ where $A$ is an open subset of $D$. Then $Pf=\chi_{PA}$ where $PA$ denotes the polarisation of the set $A$; in other words, 
\begin{equation}\label{polarization_of_set}
PA = A\cap\tau A \bigcup \big( A \cup \tau A \big)\cap(0,\,1)
\end{equation}
where $\tau:\,D\rightarrow D$ stands for the reflection $x\mapsto -x$. We shall sometimes refer to the symmetric resp. non-symmetric parts of $PA$; that is,
\begin{equation}\label{symmetric_and_nonsymetric_parts_of_PA}
\begin{array}{lcl}
A_1 & := & A\cap\tau A;\\
A_2 & := & \big(A \cup \tau A \big)\cap(0,\,1)\setminus A\cap\tau A.
\end{array}
\end{equation}

\medskip

\begin{lemma}\label{properties_of_polarised_function}
Let $f\in U_t$ for some $t\in(0,\,2/3)$. Then the following statements are equivalent:
\begin{itemize}
\item[{\it (i)}]
$f\in PU_t$;
\item[{\it (ii)}]
$f=1$ on $S:=\left\{ x\in (0,\,1):\,f(-x)=1\,\right\}$.
\end{itemize}
\end{lemma}

\noindent{\em Proof.}
Let $f\in PU_t$ so that $f=P g$ for some $g\in U_t$. Let $x\in(0,\,1)$ with $f(-x)=1$. Then $1=f(-x)=P g(-x) = g(x)\wedge g(-x)$. So $g(x)=1$ and $f(x) = P g(x) = g(x)\vee g(-x)=1$.
On the other hand, suppose that $f=1$ on $S$. For $x\in S$,
\[
P f(x)=1\vee f(-x)=1=f(x)
\text{ while }
P f(-x)=1\wedge f(-x)=f(-x),
\]
and for $x\in (0,\,1)\setminus S$,
\[
P f(x)= f(x)\vee 0=f(x)
\text{ while }
P f(-x)= f(x) \wedge 0=0=f(-x).
\]
In other words, $f=P f$. 
\qed

\medskip

\noindent It is sometimes useful to polarise with respect to the closed half-space $(-\infty,\,0]$. To distinguish between these two polarisations we use the notations $P_+$, $P_-$. In particular,
\begin{equation}\label{variant_polarisation_wrt_negative_halfspace}
P_-f(x) := 
\left\{
\begin{array}{lcr}
f(x)\wedge f(-x) & \text{ if } & 0 < x < 1,\\
f(x)\vee   f(-x) & \text{ if } & -1 < x \leq 0.\\
\end{array}
\right.
\end{equation}

\medskip

\begin{lemma}\label{polarisation_of_f_and_g}
Let $f=\chi_A\in P_+U^{(m)}_t$ for some $m\in\mathbb{N}$ and $t\in(0,\,2/3)$ where $A$ satisfies condition (\ref{condition_on_A}). Put $g:=\chi_B$ where $B:=D\setminus\overline{A}$. Then $g$ is an $m$-version of $1-f$ and $g\in P_-U_{3/2 - t}$. 
\end{lemma}

\noindent{\em Proof.}
We may suppose that $A=\bigcup_{j=1}^k A_j$ for some $1\leq k\leq m$ and $A_j=(a_j,\,b_j)$ with
\[
-1\leq a_1 < b_1 < a_2 < b_2 < \cdots < a_{k} < b_{k} \leq 1.
\]
We use the criterion in Lemma \ref{properties_of_polarised_function}. Let $x\in(-1,\,0)$ such that $g(-x)=1$. We first note that $x$ cannot be a boundary point (that is, $x\not\in\left\{a_1,\,\ldots,a_k,\,b_1,\ldots,b_k\,\right\}$). For if it is, then either $-x$ is a boundary point or $-x\in A$. This is due to the fact that $f$ is polarised to the right. In either case, we obtain the contradiction that $g(-x)=0$. We want to show that $g(x)=1$ so suppose on the contrary that $g(x)=0$. Then for $y=-x\in(0,\,1)$, it holds that $f(-y)=1$, but $f(y)=0$. This counters the fact that $f\in P_+U_t$ by the criterion. 
\qed

\section{A (non-)optimality criterion}

\noindent In this section we develop a (non-)optimality criterion for configurations $f$ in $U^{(m)}_t$. Given $f\in U_t$ define $u:=Gf$. It is known that $D(-\Delta)=W^{1,2}_0(D)\cap W^{1,2}(D)$. Thus, $u\in W^{2,2}(D)$ and by a Sobolev inequality (see \cite{Evans1998} 5.6.3 for example), $u$ belongs to the H\"{o}lder space $C^{1,\,1/2}(\overline{D})$. Define
\[
h=h_f:=\frac{u}{\psi}.
\]
Then $h\in C(D)$ and by l'H\^{o}pital's rule,
\begin{equation}\label{h_and_u_prime_at_minus_1}
h(-1) = \lim_{x\downarrow -1}\frac{u^\prime(x)}{-x}=u^\prime(-1),
\end{equation}
and similarly $h(1)=-u^\prime(1)$ at the right-hand end-point. In short, $h\in C(\overline{D})$.

\medskip

\begin{lemma}\label{l_Hopital_type_result}
Suppose that $f=\chi_A$ for some open subset $A$ in $D$. Let $u\in C(D)$. Given $a\in\overline{A}\cap D$, put $A_\eta:=[\,a-\eta,\,a+\eta]$ for $\eta>0$ small. Then
\begin{itemize}
\item[{\it (i)}]
$\lim_{\eta\downarrow 0}\frac{(f\chi_{A_\eta},\,u)}{(f\chi_{A_\eta},\,\psi)}=h(a)$;
\item[{\it (ii)}]
$\lim_{\eta\downarrow 0}\frac{(f\chi_{A_\eta},\,G[f\chi_{A_\eta}])}{(f\chi_{A_\eta},\,\psi)}=0$.
\end{itemize}
\end{lemma}

\noindent{\em Proof.}
{\it (i)} Notice that $A\cap(a-\eta,\,a+\eta)\neq\emptyset$ for each $\eta>0$. Consequently, $(f\chi_{A_\eta},\,1)=m(A\cap A_\eta)>0$ for each $\eta>0$ (small) and likewise for $(f\chi_{A_\eta},\,\psi)$. Write
\[
\frac{(f\chi_{A_\eta},\,u)}{(f\chi_{A_\eta},\,\psi)}
=
\frac{u(a)\,(f\chi_{A_\eta},\,1) + (f\chi_{A_\eta},\,u-u(a))}{
\psi(a)\,(f\chi_{A_\eta},\,1) + (f\chi_{A_\eta},\,\psi-\psi(a))}
= h(a) + \frac{\psi(a)\zeta_1 - u(a)\zeta_2}{\psi(a)(\psi(a)+\zeta_2)}
\]
where
\begin{eqnarray*}
\zeta_1 & := & \frac{(f\chi_{A_\eta},\,u-u(a))}{(f\chi_{A_\eta},\,1)},\\
\zeta_2 & := & \frac{(f\chi_{A_\eta},\,\psi-\psi(a))}{(f\chi_{A_\eta},\,1)}.
\end{eqnarray*}
Both these last vanish in the limit $\eta\downarrow 0$ and this leads to the identity.

\medskip

\noindent{\it (ii)}
From the estimate (\ref{bound_for_G_in_terms_of_psi}), for $\eta>0$ small,
\[
\psi^{-1}\,G[f\chi_{A_\eta}]\leq\psi^{-1}\,G[\chi_{A_\eta}]\leq 2\,\eta,
\]
and this establishes the limit. 
\qed

\medskip

\noindent With this preparation in hand we arrive at the crucial (non-)optimality condition. 

\begin{theorem}\label{optimality_criterion}
Let $t\in(0,\,2/3)$, $m\in\mathbb{N}$ and $f=\chi_A\in U^{(m)}_t$. Assume that $A$ satisfies condition (\ref{condition_on_A}). Suppose that $a,\,b\in D$ with $a\neq b$ such that
\begin{itemize}
\item[{\it (i)}]
$h(a) < h(b)$;
\item[{\it (ii)}]
$a\in\partial A$;
\item[{\it (iii)}]
$b\in\partial A$.
\end{itemize}
Then there exists $f_1\in U^{(m)}_t$ such that
$J(f) < J(f_1)$.
\end{theorem}

\noindent{\em Proof.}
Write $f$ in the form $f=\sum_{j=1}^k\chi_{A_k}$
for some $1\leq k\leq m$ where $A_j=(a_j,\,b_j)$ and
\[
-1\leq a_1 < b_1 < a_2 < b_2 < \cdots < a_k < b_k \leq 1.
\]
Given $\eta>0$ put $A_\eta:=[a-\eta,\,a+\eta]$ and $B_\eta:=(b-\eta,\,b+\eta)$. Set $g:= 1-f$. The functions
\[
\eta\mapsto(\chi_{A_\eta}\,f,\,\psi)\text{ and }\eta\mapsto(\chi_{B_\zeta}\,g,\,\psi)
\]
are strictly increasing at least for $\eta>0$ small. For $\varepsilon>0$ sufficiently small, there exist unique $\eta>0$ and $\zeta>0$ depending upon  $\varepsilon$ such that
\[
\varepsilon = (\chi_{A_\eta}\,f,\,\psi) =  (\chi_{B_\zeta}\,g,\,\psi).
\]
Define
\[
f_\varepsilon:= f - f\,\chi_{A_\eta} + g\,\chi_{B_\zeta}.
\]
Then $f_\varepsilon\in U^{(m)}_t$ for $\varepsilon>0$ small. Now
\begin{eqnarray}\label{computation_for_J_of_f_epsilon_minus_J_of_f}
\nonumber
J(f_\varepsilon) - J(f) & = & (f_\varepsilon,\,G\,f_\varepsilon)-(f,\,G\,f)\\ \nonumber
                        & = & (f_\varepsilon - f,\,G\,\big[\,f_\varepsilon+f\,\big])\\ \nonumber
                        & = &    (-f\chi_{A_\eta}+g\chi_{B_\zeta},\,G\,\big[\,2\,f-f\chi_{A_\eta}+ g\chi_{B_\zeta}\,\big])\\ 
                        & = &  2\,(g\chi_{B_\zeta},\,u) - 2\,(f\chi_{A_\eta},\,u)+ (g\chi_{B_\zeta}-f\chi_{A_\eta},\,G\,\big[\,g\chi_{B_\zeta}-f\chi_{A_\eta}\,\big])
\end{eqnarray}
where $u=G\,f$ as usual. Thus, by Lemma \ref{l_Hopital_type_result} (with the help of the Cauchy-Schwarz inequality to deal with the cross-terms),
\[
\lim_{\varepsilon\downarrow 0}\varepsilon^{-1}\Big\{\,J(f_\varepsilon) - J(f)\Big\}=2(\,h(b) - h(a)\,)>0.
\]
In particular, there exists $\varepsilon>0$ (small) such that $f_1:=f_\varepsilon\in U^{(m)}_t$ satisfies $J(f_1)>J(f)$.
\qed

\section{More on the (non-)optimality condition}
 
\noindent In this section we verify condition {\it (i)} in Theorem \ref{optimality_criterion} for some particular configurations $f$ in $U_t$. 

\begin{lemma}\label{evaluation_of_u_prime}
Let $f\in U_t$ for some $t\in(0,\,2/3)$ and set $u:=G\,f$. Then
\begin{itemize}
\item[{\it (i)}]
$u^\prime(-1) + u^\prime(1) = -\int_D x\,f\,dm$; 
\item[{\it (ii)}]
$u^\prime(-1) - u^\prime(1) = \int_D f\,dm$.
\end{itemize}
\end{lemma}

\medskip

\noindent{\em Proof.}
{\it (i)} Using $\psi^\prime=-x$ and the integration-by-parts formula, 
\begin{eqnarray*}
\int_D x\,f\,dm & = & \int_D\psi^\prime\,u^{\prime\prime}\,dm\\
 & = & \psi^\prime(1)\,u^\prime(1) - \psi^\prime(-1)\,u^\prime(-1)
       - \int_D\psi^{\prime\prime}\,u^{\prime}\,dm\\
 & = & - u^\prime(1) - u^\prime(-1).      
\end{eqnarray*}
{\it (ii)} This follows from
$
\int_D u^{\prime\prime}\,dm = -\int_D f\,dm = u^\prime(1) - u^\prime(-1).
$
\qed

\medskip

\noindent  Assuming that $f\not\equiv 0$, define $a,\,b\in\mathbb{R}$ by
\begin{equation}\label{definition_of_a_and_b}
\begin{array}{lcl}
a     & := & \inf\Big\{\,x>-1:\,(\,\chi_{(-1,\,x]},\,f)>0\,\Big\},\\
1 - b & := & \inf\Big\{\,x>0:\,(\,\chi_{[1-x,\,1)},\,f)>0\,\Big\};
\end{array}
\end{equation}
so that $a\in[-1,\,1)$ and $b\in(-1,\,1]$. 

\medskip

\begin{lemma}\label{evaluation_of_h}
For each $f\in U_t$ with $f\not\equiv 0$,
\begin{itemize}
\item[{\it (i)}]
$h(y) = \frac{1}{1-y}\int_D\big(\,1 - x\,\big)\,f\,dm$ for each $y\in[-1,\,a]$; 
\item[{\it (ii)}]
$h(y) = \frac{1}{1+y}\int_D\big(\,1 + x\,\big)\,f\,dm$ for each $y\in[b,\,1]$.
\end{itemize}
\end{lemma}

\medskip

\noindent{\em Proof.}
{\it (i)} Suppose that $a=-1$. By (\ref{h_and_u_prime_at_minus_1}) and Lemma \ref{evaluation_of_u_prime}, 
\[
h(-1) = u^\prime(-1)=(1/2)\,\int_D\big(\,1 - x\,\big)\,f\,dm.
\]
Now suppose that $a\in(-1,\,1)$. Using integration-by-parts,
\begin{eqnarray*}
u(y) & = & \int_{(-1,\,y]}u^\prime\,dm\\
     & = & -\int_{(-1,\,y]}u^\prime\,\psi^{\prime\prime}dm\\
     & = & -u^\prime(y)\,\psi^\prime(y) + u^\prime(-1)            \,\psi^\prime(-1) + \int_{(-1,\,y]}u^{\prime\prime}\,\psi^{\prime}dm\\
     & = & -u^\prime(y)\,\psi^\prime(y) + u^\prime(-1)           
\end{eqnarray*}
as $u^{\prime\prime}=-f=0$ $m$-a.e. on $(-1,\,a]$. For the same reason,
\[
u^\prime(y) - u^\prime(-1) = \int_{(-1,\,y]}u^{\prime\prime}\,dm = 0.
\]
Therefore,
\[
u(y) = \big(\,1 - \psi^\prime(y)\,\big)\,u^\prime(-1)
     = \frac{1+y}{2}\,\int_D\big(\,1 - x\,\big)\,f\,dm
\]
from which the statement is clear. Part {\it (ii)} follows in a similar fashion. 
\qed

\medskip

\begin{proposition}\label{inequality_between_h_at_a_and_h_at_b}
Let $f\in P U_t$ for some $t\in(0,\,2/3)$. With $a,\,b$ as in (\ref{definition_of_a_and_b}) assume that
\begin{itemize}
\item[{\it (i)}]
$-1 < a < 0 < -a < b < 1$;
\item[{\it (ii)}]
$f=1$ $m$-a.e. on $(-a,\,b)$.
\end{itemize}
Then $h(a) < h(b)$. 
\end{proposition}

\noindent{\em Proof.}
From Lemma \ref{evaluation_of_h},
\[
h(b) - h(a) = -\frac{a+b}{(1-a)\,(1+b)}\int_D f\,dm
+ \frac{2-a+b}{(1-a)\,(1+b)}\int_D x\,f\,dm.
\]
Put $S:=\left\{ x\in (0,\,1):\,f(-x)=1\,\right\}$ as in Lemma \ref{properties_of_polarised_function}. Then
\[
\int_D x\,f\,dm = \int_S x\left\{\,f(x) - f(-x)\,\right\}\,m(dx)
+
\int_{(0,\,1)\setminus S} x\,f\,dm
\]
\[
=
\int_{(0,\,1)\setminus S} x\,f\,dm
\geq  
\int_{(-a,\,b)} x\,f\,dm
=(1/2)( b^2 - a^2 )
\]
making use of {\it (ii)}. Thus,
\[
(2-a+b)\,\int_D x\,f\,dm \geq (1/2)\,(2-a+b)\,( b^2 - a^2 ) > b^2 - a^2 \geq (a+b)\,\int_D f\,dm
\]
and hence $h(b) - h(a)>0$. 
\qed

\medskip

\noindent In the next two sections, we show non-optimality of polarised configurations in three broad cases. 
\section{Two non-symmetric cases}

\noindent Let $t\in(0,\,1/3]$ and imagine a configuration polarised to the right that charges the left-hand interval $(-1,\,0)$ but which is not symmetric under reflection in the origin. We show this is non-optimal. 

\medskip

\begin{lemma}\label{lemma_on_non_symmetric_case_I}
Let $m\in\mathbb{N}$ and $t\in(0,\,1/3]$. Suppose that $f\in U^{(m)}_t$ satisfies the properties
\begin{itemize}
\item[{\it (i)}]
$f=Pf$;
\item[{\it (ii)}]
$(\,f,\,\chi_{(-1,\,0)}\,)>0$.
\item[{\it (iii)}]
$(\,f,\,\chi_{(0,\,1)}\,)>(\,f,\,\chi_{(-1,\,0)}\,)$.
\end{itemize}
Then there exists $f_1\in U^{(m)}_t$ with the property that $f_1 = P\,f_1$ such that $J(f) < J(f_1)$.
\end{lemma}

\noindent{\em Proof.}
We may assume that $f=\chi_A$ where $A$ satisfies condition (\ref{condition_on_A}). We may then write $f$ in the form described at the beginning of the proof of Theorem \ref{optimality_criterion}. By {\it (ii)}, $a_1<0$; and by {\it (i)}, $b_k\geq-a_1$.

\medskip

\noindent{\em Case (a):} $-1 < a_1$ and $b_k<1$. Then, in fact,
$-1 < a_1 < 0 < -a_1 \leq b_k < 1$. Put 
\[
k_1:=\min\Big\{\,j = 1,\ldots, k:\,b_j \geq -a_1\,\Big\}.
\]
Suppose first of all that $-a_1=b_{k_1}$. Decompose $A$ into its symmetric and non-symmetric parts $A_1$ and $A_2$ as in (\ref{symmetric_and_nonsymetric_parts_of_PA}). By {\it (iii)}, $A_2\neq\emptyset$. Write $f_1:=\chi_{A_1}$ and $f_2:=\chi_{A_2}$. By symmetry, $h_{f_1}(a_1)=h_{f_1}(-a_1)$. Further, $h_{f_2}(a_1)<h_{f_2}(-a_1)$, this being a consequence of (\ref{Green_kernel_inequality}). Therefore, as $h=h_{f_1}+h_{f_2}$, we obtain $h(a_1)<h(-a_1)$. The conclusion follows with an application of Theorem \ref{optimality_criterion}. 

\medskip

\noindent If $-a_1\neq b_{k_1}$ then $-a_1<b_{k_1}$ and $f=1$ on $(-a_1,\,b_{k_1})$. If $k_1 =k$ then $h(a_1) < h(b_k)$ by Proposition \ref{inequality_between_h_at_a_and_h_at_b}. On the other hand, if $k_1 < k$ define
\[
f_1:=\sum_{j=1}^{k_1}\chi_{A_j}\text{ and }f_2:=\sum_{j=k_1+1}^{k}\chi_{A_j}.
\]
By Proposition \ref{inequality_between_h_at_a_and_h_at_b},
$h_{f_1}(a_1) < h_{f_1}(b_{k_1})$. It can be seen from the representation in Lemma \ref{evaluation_of_h} that $h_{f_2}$ is increasing on $[-1,\,a_{k_1+1}]$. In sum, then, $h(a_1) < h(b_{k_1})$. Now apply Theorem \ref{optimality_criterion} once more. 

\medskip

\noindent{\em Case (b):} $-1 < a_1$ and $b_k = 1$. In this situation,
$-1 < a_1 < 0 < -a_1 < b_k = 1$. Define $k_1$ as before. The case $k_1<k$ may be dealt with in a similar way to case (a) above. So assume that $k_1=k$. 
As $f$ is polarised to the right, the interval $((-b_1)\vee 0,\,-a_1)$ must sit inside $A$ and so it must hold that $a_k<-a_1$. In case $t\in(0,\,1/3)$, it must also hold that $0 < a_k$. The situation $t=1/3$ and $a_k=0$ forces $(\,f,\,\chi_{(-1,\,0)}\,)=0$ contradicting {\it (ii)}. In either case, therefore, $0< a_k < -a_1$ and $k\geq 2$. 

\medskip

\noindent Consider the function $g:=\chi_B$ where $B:=D\setminus\overline{A}$. By Lemma \ref{polarisation_of_f_and_g}, $g\in P_-U_{3/2-t}$. Thus,
\[
-1 < a_1 < - a_k < 0 < a_k < b_k = 1,
\]
and $g=1$ just to the right of $-a_k$ as $g$ is polarised to the left. This situation corresponds to the one described at the start of the consideration of this case but for $g$ instead of $f$. Use the fact that $h_g = 1 - h_f$.

\medskip

\noindent{\em Case (c):} $a_1 = -1$. Then $b_k=1$ and $a_k\leq -b_1$ as $f=Pf$. Apply the arguments in case {\em (a)} to the function $g$.  
\qed

\medskip

\noindent We now take $t\in(0,\,1/3)$ and imagine a configuration that lies entirely in the right-hand interval $(0,\,1)$ but that has not yet been pushed rightwards to the maximum extent. We again show non-optimality. 

\medskip

\begin{lemma}\label{lemma_on_non_symmetric_case_II}
Let $m\in\mathbb{N}$ and $t\in(0,\,1/3)$. Suppose that $f\in U^{(m)}_t$ satisfies the properties
\begin{itemize}
\item[{\it (i)}]
$(\,f,\,\chi_{(-1,\,0)}\,)=0$;
\item[{\it (ii)}]
$(\,f,\,\chi_{(0,\,\xi_t)}\,)>0$.
\end{itemize}
Then there exists $f_1\in U^{(m)}_t$ with the property that $f_1 = P\,f_1$ such that $J(f) < J(f_1)$.
\end{lemma}

\noindent{\em Proof.}
Again take $f=\chi_A$ where $A$ satisfies condition (\ref{condition_on_A}) and
suppose $f$ takes the form described at the beginning of the proof of Theorem \ref{optimality_criterion}. By {\it (i)}, $a_1\geq 0$ and by {\it (ii)}, $a_1<\xi_t$.
Therefore $a_1<b_1<1$; for otherwise, if $b_1=1$ then 
\[
(f,\,\psi) \geq (\chi_{(a_1,\,1)},\,\psi) > (\chi_{(\xi_t,\,1)},\,\psi) = t.
\]
Again borrowing the notation of Theorem \ref{optimality_criterion}, put
$f_1:=\chi_{A_1}$ and $f_2:=\sum_{j=2}^{k}\chi_{A_j}$. By Lemma \ref{evaluation_of_h},
\[
h_{f_1}(a_1)=\frac{1}{1-a_1}\int_{(a_1,\,b_1)}\big(\,1-x\,\big)\,dm = \frac{b_1 - a_1}{1-a_1}\left\{\,1 - (1/2)\big( a_1 + b_1 \big)\right\}
\]
and
\[
h_{f_1}(b_1)=\frac{1}{1+b_1}\int_{(a_1,\,b_1)}\big(\,1+x\,\big)\,dm = \frac{b_1 - a_1}{1+b_1}\left\{\,1 + (1/2)\big( a_1 + b_1 \big)\right\}.
\]
A little algebra yields $h_{f_1}(b_1) > h_{f_1}(a_1)$. Lemma \ref{evaluation_of_h} also indicates that $h_{f_2}$ is monotone increasing on $[-1,\,a_2]$. Therefore, $h(b_1) > h(a_1)$.
The conclusion now follows with the help of Theorem \ref{optimality_criterion}.
\qed

\section{The symmetric case}

\noindent In the last of the three cases, we consider a configuration that is symmetric under reflection in the origin. 

\begin{proposition}\label{symmetric_case}
Let $m\in\mathbb{N}$ and $t\in(0,\,2/3)$. Suppose that $f\in U^{(m)}_t$ satisfies the properties
\begin{itemize}
\item[{\it (i)}]
$f = Pf$;
\item[{\it (ii)}]
$(\,f,\,\chi_{(-1,\,0)}\,)>0$;
\item[{\it (iii)}]
$(\,f,\,\chi_{(0,\,1)}\,)=(\,f,\,\chi_{(-1,\,0)}\,)$.
\end{itemize}
Then there exists $f_1\in U^{(m)}_t$ with the property that $f_1 = P\,f_1$ such that $J(f) < J(f_1)$.
\end{proposition}

\medskip

\noindent Before embarking on the proof of Proposition \ref{symmetric_case}, we require a number of supplementary results. 

\medskip

\begin{lemma}\label{on_diaginal_matrix_element}
Suppose that $f=\chi_A\in U^{(m)}_t$ for some $m\in\mathbb{N}$ and $t\in(0,\,2/3)$. Assume that $A$ satisfies condition (\ref{condition_on_A}). Suppose that $a\in\partial A\cap D$. Then
\[
\lim_{\eta\downarrow 0}\frac{(f\chi_{A_\eta},\,G[f\chi_{A_\eta}])}{(f\chi_{A_\eta},\,\psi)^2}=\psi(a)^{-1},
\]
where $A_\eta=[a-\eta,\,a+\eta]$ as before.
\end{lemma}

\noindent{\em Proof.}
Write
\[
G[f\chi_{A_\eta}](x)=
\psi(a)\,(f\chi_{A_\eta},\,1) +
\left\{\,\psi(x) - \psi(a)\,\right\}(f\chi_{A_\eta},\,1)
+ r(x)
\]
where $r(x):=(\,G(x,\cdot) - \psi(x),\,f\chi_{A_\eta}\,)$
for $x\in D$. Since $\psi(x)=G(x,\,x)$, the estimate (\ref{variation_estimate_for_Green_function}) gives
\[
\big|\,(\,f\chi_{A_\eta},\,r)\big| \leq 2\,\eta\,(f\chi_{A_\eta},\,1)^2.
\]
Forming the inner product we obtain
\[
(f\chi_{A_\eta},\,G[f\chi_{A_\eta}])
=
\psi(a)\,(f\chi_{A_\eta},\,1)^2 + (f\chi_{A_\eta},\,\psi - \psi(a))\,(f\chi_{A_\eta},\,1) + (\,f\chi_{A_\eta},\,r).
\]
It is clear from this that 
\[
\lim_{\eta\downarrow 0}
\frac{(f\chi_{A_\eta},\,G[f\chi_{A_\eta}])}{(f\chi_{A_\eta},\,1)^2}=\psi(a).
\]
A short step leads to the assertion. 
\qed

\begin{lemma}\label{off_diaginal_matrix_element}
Suppose that $f=\chi_A\in U^{(m)}_t$ for some $m\in\mathbb{N}$ and $t\in(0,\,2/3)$ and that $A$ satisfies condition (\ref{condition_on_A}). Suppose that $a,\,b\in D$ with $a\neq b$ such that both $a\in\partial A$ and $b\in\partial A$. Put $g:=1-f$. Given $\varepsilon>0$ sufficiently small
there exist unique $\eta>0$ and $\zeta>0$ depending upon $\varepsilon$ such that
$\varepsilon = (f\,\chi_{A_\eta},\,\psi) =  (g\,\chi_{B_\zeta},\,\psi)$. Then
\[
\lim_{\varepsilon\downarrow 0}\frac{(f\chi_{A_\eta},\,G[g\chi_{B_\zeta}])}{\varepsilon^2}=\frac{G(a,\,b)}{\psi(a)\,\psi(b)}.
\]
\end{lemma}

\noindent{\em Proof.}
Write
\[
G[g\chi_{B_\zeta}](x) = (G(x,\cdot),\,g\chi_{B_\zeta})
=G(a,\,b)\,(g\chi_{B_\zeta},\,1) + r(x)
\]
where $r(x):=(G(x,\cdot)-G(a,\,b),\,g\chi_{B_\zeta})$.
For $x\in A_\eta$ and $y\in B_\zeta$,
\[
\big|\,G(x,\,y) - G(a,\,b)\,\big| \leq \eta + \zeta
\]
by (\ref{variation_estimate_for_Green_function}). Consequently,
\[
\big|\,(\,f\chi_{A_\eta},\,r)\big| \leq \big(\,\eta + \zeta\,\big)\,(f\chi_{A_\eta},\,1)\,(g\chi_{B_\zeta},\,1).
\]
Now,
\[
(f\chi_{A_\eta},\,G[g\chi_{B_\zeta}]) =G(a,\,b)\,(f\chi_{A_\eta},\,1)\,(g\chi_{B_\zeta},\,1)
+(\,f\chi_{A_\eta},\,r),
\]
from which we derive
\[
\lim_{\varepsilon\downarrow 0}
\frac{(f\chi_{A_\eta},\,G[g\chi_{B_\zeta}])}{(f\chi_{A_\eta},\,1)\,(g\chi_{B_\zeta},\,1)}=G(a,\,b),
\]
and the conclusion follows straightforwardly. 

\qed

\begin{lemma}\label{technical_lemma_1_for_symmetry_case}
Let $b\in(0,\,1)$ and $a:=-b$. Let $\eta>0$ small and define $\zeta=\zeta(\eta)$ via the relation
\[
(\chi_{[a,\,a+\eta]},\,\psi) = (\chi_{[b,\,b+\zeta]},\,\psi).
\] 
Then $\zeta$ depends smoothly upon $\eta$ in a neighbourhood of $\eta=0$ and
\[
\zeta = \eta + \frac{2\,b}{1-b^2}\,\eta^2 +O(\eta^3)
\]
in the limit $\eta\downarrow 0$. 
\end{lemma}

\noindent{\em Proof.}
A short computation gives that
\[
(\chi_{[a,\,a+\eta]},\,\psi) = (1/2)\left\{\,\big( 1 - a^2 \big)\,\eta - a\,\eta^2 - (1/3)\,\eta^3\,\right\}.
\]
Define smooth functions $f,\,g:\mathbb{R}\rightarrow\mathbb{R}$ by
\[
f(\eta)  := \big( 1 - b^2 \big)\,\eta + b\,\eta^2 - (1/3)\,\eta^3
\text{ and }
g(\zeta) := \big( 1 - b^2 \big)\,\zeta - b\,\zeta^2 - (1/3)\,\zeta^3.
\]
Now $f^\prime(0) = g^\prime(0) = 1-b^2>0$. In particular, $f$ is strictly increasing in a neighbourhood of $\eta = 0$ and $g$ possesses a local smooth inverse $h$ in the neighbourhood of $\zeta = 0$ by the inverse function theorem. Note that $\zeta$ is characterised by the relation $g(\zeta)=f(\eta)$ for $\eta>0$ small. Thus $\zeta = (h\circ f)(\eta)$ and depends smoothly upon $\eta$ for $\eta>0$ small. Implicit differentiation yields $\zeta^\prime(0)=1$ and $\zeta^{\prime\prime}(0)=\frac{4\,b}{1-b^2}$. Taylor's theorem with remainder then yields the expansion. 
\qed

\medskip

\noindent{\em Proof of Proposition \ref{symmetric_case}.} 
We may suppose that $f=\chi_A$ where $A$ satisfies condition (\ref{condition_on_A}). Define $a$ as in (\ref{definition_of_a_and_b}). Assume in the first instance that $a\in(-1,\,0)$. Put $b:=-a$. Conditions {\it (i)}-{\it (iii)} entail that $f$ is even. In particular, $a,\,b\in\partial A$. As in Lemma \ref{off_diaginal_matrix_element} we write
\[
\varepsilon = (f\,\chi_{A_\eta},\,\psi) = (g\,\chi_{B_\zeta},\,\psi)
\]
for $\varepsilon>0$ small. We aim to show that $J(f_\varepsilon) - J(f)>0$ at least for $\varepsilon>0$ small as in Theorem \ref{optimality_criterion} and we shall borrow notation without comment from its proof. We first claim that (see (\ref{computation_for_J_of_f_epsilon_minus_J_of_f}))
\begin{eqnarray}\label{claim_in_proof_of_symmetry_theorem}
\nonumber
\lim_{\varepsilon\downarrow 0}\frac{J(f_\varepsilon) - J(f)}{\varepsilon^2}
& = & 
\lim_{\varepsilon\downarrow 0}\,\varepsilon^{-2}
\left\{
2\,(g\chi_{B_\zeta},\,u) - 2\,(f\chi_{A_\eta},\,u)+ (g\chi_{B_\zeta}-f\chi_{A_\eta},\,G\,\big[\,g\chi_{B_\zeta}-f\chi_{A_\eta}\,\big])
\right\}\\
& = & 
\frac{2}{\psi(b)^2}\left\{\,
b\,h(b)
+ u^\prime(b)
+ b\,(1-b)\,\right\}.
\end{eqnarray}

\medskip

\noindent As $u\in C^{1,\,1/2}(\overline{D})$ we have that
\[
( f\chi_{A_\eta},\,u ) = u(a)\,\eta + u^\prime(a)\,(1/2)\eta^2 + O(\eta^{5/2}).
\]
Also, with the help of Lemma \ref{technical_lemma_1_for_symmetry_case},
\begin{eqnarray*}
(g\,\chi_{B_\zeta},\,u) & = & u(b)\,\zeta + u^\prime(b)\,(1/2)\,\zeta^2 + O(\zeta^{5/2})\\
& = & u(b)\,\left\{\,\eta + \frac{2\,b}{1-b^2}\,\eta^2\,\right\}
+ u^\prime(b)\,(1/2)\,\eta^2 + O(\eta^{5/2})\\
& = & u(b)\,\eta + 
\left\{\,\frac{2\,b}{1-b^2}\,u(b)
+ \frac{1}{2}u^\prime(b)\,\right\}
\,\eta^2 + O(\eta^{5/2}).
\end{eqnarray*}
Now $u^\prime(a)=-u^\prime(b)$ because $u$ is even so
\[
\lim_{\varepsilon\downarrow 0}\frac{(g\,\chi_{B_\zeta},\,u) - (f\,\chi_{A_\eta},\,u) }{\varepsilon^2}
=
\lim_{\varepsilon\downarrow 0}\big(\frac{\eta}{\varepsilon}\big)^2\frac{(g\,\chi_{B_\zeta},\,u) - (f\,\chi_{A_\eta},\,u) }{\eta^2}
=
\frac{1}{\psi(b)^2}
\left\{
b\,h(b)
+ u^\prime(b)
\right\}.
\]
On the other hand, from Lemmas \ref{on_diaginal_matrix_element} and \ref{off_diaginal_matrix_element},
\[
\lim_{\varepsilon\downarrow 0}
\frac{
(g\chi_{B_\zeta}-f\chi_{A_\eta},\,G\,\big[\,g\chi_{B_\zeta}-f\chi_{A_\eta}\,\big])
}{\varepsilon^2} = \frac{2}{\psi(a)} - \frac{2\,G(a,\,b)}{\psi(a)\,\psi(b)}
= \frac{2\,b\,(1-b)}{\psi(b)^2}
\]
as $\psi(b) - G(a,\,b) = b\,(1-b)$. The combination of these identities establishes the claim (\ref{claim_in_proof_of_symmetry_theorem}). 

\medskip

\noindent We now show that the expression in (\ref{claim_in_proof_of_symmetry_theorem}) is positive. From Lemma \ref{evaluation_of_h} and the even property of $f$, 
\[
h(b) = \frac{1}{1+b}\int_D f\,dm
\]
and from Lemma \ref{evaluation_of_u_prime},
\[
u^\prime(b) = u^\prime(1) = -\frac{1}{2}\int_D f\,dm.
\]
So we may write
\begin{eqnarray*}
b\,h(b)
+ u^\prime(b)
+ b\,(1-b)
& = &
\left\{\,\frac{b}{1+b}-\frac{1}{2}\,\right\}\int_D f\,dm
+ b\,(1-b)\\
& \geq & 
\left\{\,\frac{b}{1+b}-\frac{1}{2}\,\right\}2b
+ b\,(1-b)\\
& = & \frac{b^2(1-b)}{1+b}
>0.
\end{eqnarray*}
The conclusion now follows by Theorem \ref{optimality_criterion}.

\medskip

\noindent The case $a=-1$ may be dealt with by applying the above argument (with appropriate modifications) to $g:=1-f$.
\qed

\section{The main result}

\noindent We are now in a position to prove the main result Theorem \ref{first_main_theorem}. 

\begin{theorem}\label{optimiser_for_beta_m_t}
Let $t\in (0,\,1/3]$ and $m\in\mathbb{N}$ and put $A_t:=(\xi_t,\,1)$ with $\xi_t$ as in (\ref{definition_of_xi_t}). Then $\alpha^{(m)}_t = J(f)$ where $f=\chi_{A_t}$. 
\end{theorem}

\noindent{\em Proof.}
By Theorem \ref{existence_of_maximiser_for_J_in_U_m_t} there exists $f\in U^{(m)}_t$ such that $\beta^{(m)}_t = J(f)$. Now $U^{(m)}_t$ is closed under polarisation. So $P\,f\in U^{(m)}_t$ and $J(f) \leq J(P f)$ by Theorem \ref{polarisation_theorem}. We may assume therefore that $f=P f$. 

\medskip 

\noindent Assume that $(f,\,\chi_{(-1,\,0)})>0$. Since $f=Pf$ it must be the case that $(f,\,\chi_{(0,\,1)})\geq(f,\,\chi_{(-1,\,0)})$ (in consequence of Lemma \ref{properties_of_polarised_function}). By Lemma \ref{lemma_on_non_symmetric_case_I} and Proposition \ref{symmetric_case}, there exists $f_1\in U^{(m)}_t$ with the property that $f_1 = P f_1$ and $J(f) < J(f_1)$. This contradicts the optimality of $f$. We conclude that $(f,\,\chi_{(-1,\,0)})=0$. 

\medskip

\noindent If $t=1/3$ this compels $f=\chi_{(0,\,1)}$ bearing in mind that $f\in U^{(m)}_{1/3}$ and $f=P\,f$. So let us now take $t\in (0,\,1/3)$. Suppose that $(f,\,\chi_{(0,\,\xi_t)})>0$. Then the requirements of Lemma \ref{lemma_on_non_symmetric_case_II} are satisfied and hence there exists $f_1\in U^{(m)}_t$ with the property that $f_1 = P\,f_1$ such that $J(f) < J(f_1)$. Again this contradicts optimality. Hence, $(f,\,\chi_{(0,\,\xi_t)})=0$. In fact, $(f,\,\chi_{(-1,\,\xi_t)})=0$. As $f\in U^{(m)}_t$ we draw the conclusion that $f=\chi_{A_t}$. 
\qed

\medskip

\begin{corollary}\label{beta_t_for_t_in_0_one_third}
Let $t\in (0,\,1/3]$. Then $\alpha_t = J(f)$ where $f=\chi_{A_t}$. 
\end{corollary}

\medskip

\noindent{\em Proof.}
Let $f\in U_t$. By Lindel\"{o}f's theorem, we may write $f$ in the form $f=\chi_A$ where $A=\bigcup_{k=1}^\infty A_k$ is a countable union of disjoint open intervals $A_k$ in $D$. Put $f_n:=\sum_{k=1}^n\chi_{A_k}$. By the monotone convergence theorem, $J(f)=\lim_{n\rightarrow\infty}J(f_n)$. Note that $f_n\in U^{(n)}_t$. By Theorem \ref{optimiser_for_beta_m_t}, $J(f_n)\leq J(\chi_{A_t})$ and $J(f)\leq J(\chi_{A_t})$ on taking limits.
\qed

\medskip

\noindent{\em Proof of Theorem \ref{first_main_theorem}.} We only need to deal with the case $t\in(1/3,\,2/3)$ in view of Corollary \ref{beta_t_for_t_in_0_one_third}. Let $f\in V_t$ for such a $t$. Then $g:=1-f\in V_{2/3 - t}$ and
\[
J(f) = 2(t - 1/3) + J(g) \leq
2(t - 1/3) + J(\chi_{A_{2/3 - t}})
= J(1 - \chi_{A_{2/3 - t}}) = J(\chi_{A_t}).
\] 
This clinches the result in the final case. \qed
\section{Application: maximum flux exchange flow}

\noindent In this section we prove Theorem \ref{second_main_theorem}. 

\medskip

\begin{proposition}
It holds that
\begin{itemize}
\item[{\it (i)}]
$\gamma_\lambda = 2\,\alpha_{\frac{1-\lambda}{3}}-(1/3)(1-\lambda)^2$
for $\lambda\in(-1,\,1)$;
\item[{\it (ii)}]
$\gamma = \sup_{\lambda\in(-1,\,1)}\gamma_\lambda$. 
\end{itemize}
\end{proposition}

\noindent{\em Proof.}
{\it (i)}. Fix $\lambda\in(-1,\,1)$. 
Let $A$ be an open subset in $D$. Suppose that $u$ satisfies (\ref{equations_for_u}) along with the flux-balance condition $(u,\,1)=0$. Put
\[
f=f_{A,\lambda}:=\left\{
\begin{array}{lll}
-(\lambda + 1) & \text{ on } & A,\\
-(\lambda - 1) & \text{ on } & D\setminus A,
\end{array}
\right.
\]  
Then $u=Gf$. From the flux-balance condition and symmetry of the Green operator,
\[
0=(1,\,G\,f)=(\psi,\,f)=-(\lambda+1)(\psi,\,\chi_{A})-(\lambda-1)(\psi,\,\chi_{D\setminus A});
\]
so that $\lambda= (\psi,1)^{-1}(\psi,\chi_{D\setminus A}-\chi_{A})$ and
$(\psi,\,\chi_A) = \frac{1-\lambda}{2}(\psi,\,1)$. 
Moreover, 
\begin{eqnarray*}
(\chi_{D\setminus A},\, u) & = & (G\,\chi_{D\setminus A},\,f)\\
             & = & -(\lambda+1)(G\,\chi_{D\setminus A},\,\chi_{A})-(\lambda-1)(G\,\chi_{D\setminus A},\,\chi_{D\setminus A})\\
             & = & 2(\psi,1)^{-1}            \Big\{-(\psi,\,\chi_{D\setminus A})(G\,\chi_{D\setminus A},\,\chi_{A})+(\psi,\,\chi_{A})(G\,\chi_{D\setminus A},\,\chi_{D\setminus A})\Big\}\\              & = & 2(\psi,\,1)^{-1}\Big\{(\psi,\,1)(G\,\chi_A,\,\chi_A)-(\psi,\,\chi_{A})^2\Big\}\\
& = & 2\,J(\chi_A) - (1/2)(1-\lambda)^2(\psi,\,1),       
\end{eqnarray*}
upon writing $\chi_{D\setminus A}=\chi_D - \chi_A$ on each occurrence in the penultimate line. Now simplify using (\ref{inner_product_of_psi_with_1}). These considerations lead to the reformulation {\it (i)}. The statement in {\it (ii)} then follows immediately.  
\qed

\medskip

\noindent{\em Proof of Theorem \ref{second_main_theorem}}.
Part {\it (i)} follows from Theorem \ref{first_main_theorem}. For $t\in(0,\,2/3)$,
\[
\alpha_t = 2(t-1/3) + \alpha_{2/3 - t},
\]
as can be seen from the proof of Theorem \ref{first_main_theorem}. Therefore, for $\lambda\in(-1,\,1)$,
\[
\gamma_{-\lambda}
=2\,\alpha_{\frac{1+\lambda}{3}}-(1/3)(1+\lambda)^2
=2\,\alpha_{\frac{1-\lambda}{3}}-(1/3)(1-\lambda)^2
=\gamma_\lambda.
\]
We now show that $\gamma_\lambda<\gamma_0$ for each $\lambda\in(0,\,1)$. 

\medskip

\noindent Let $\xi\in(-1,\,1)$ and $u:=G\,f$ where $f:=\chi_{(\xi,\,1)}$. Then
\[
u(x) =
\left\{
\begin{array}{lcr}
\frac{1}{4}(\xi - 1)^2(x+1) & \text{ if } & -1 < x \leq \xi,\\
-\frac{1}{2}x^2 + \frac{1}{4}(\xi + 1)^2\,x
+ \frac{1}{2}
-\frac{1}{4}(\xi + 1)^2
& \text{ if } & \xi \leq x < 1.
\end{array}
\right.
\]
A computation leads to
\[
J(f) = -\frac{1}{8}\,\xi^4 + \frac{1}{6}\,\xi^3 + \frac{1}{4}\,\xi^2
- \frac{1}{2}\,\xi + \frac{5}{24}.
\]
Also (see (\ref{definition_of_xi_t})), 
\[
\frac{1-\lambda}{3} = (\psi,\,f) = \varphi(\xi) = \frac{1}{6}\left\{\, 2 - 3\,\xi + \xi^3\,\right\}.
\]
Therefore,
\begin{eqnarray*}
\gamma_\lambda & = & 2\,\alpha_{\frac{1-\lambda}{3}} - (1/3)(1-\lambda)^2\\
               & = & 2\,J(f) - 3\,\varphi(\xi)^2\\
               & = &
\frac{1}{12} - \frac{1}{4}\,\xi^2 + \frac{1}{4}\,\xi^4 - \frac{1}{12}\,\xi^6\\                    & =: & h(\xi).
\end{eqnarray*}
Now $h(0)=\frac{1}{12}$ and $h(1) = 0$ and $h^\prime(\xi) =-(1/2)\,\xi\,(1-\xi^2)^2 < 0$ for $\xi\in(0,\,1)$. This shows that  $\gamma_\lambda<\gamma_0$ for each $\lambda\in(0,\,1)$ as desired. The result follows from this and {\it (i)} of the Theorem.
\qed

\bigskip

\noindent{\em Acknowledgements.} I am grateful to Professor Richard Kerswell and Professor Michiel van den Berg for several helpful discussions on this topic.

\end{document}